\RequirePackage{ifpdf}
\ifpdf % We are running pdfTeX in pdf mode
\documentclass[pdftex]{sigma}
\else
\documentclass{sigma}
\fi

\newcommand{\intprod}{\;\rule{5pt}{.3pt}\rule{.3pt}{7pt}\;}
\newcommand{\Dirac}{\mbox{$\not\!\!D$}}

\begin{document}

\allowdisplaybreaks

\renewcommand{\PaperNumber}{084}

\FirstPageHeading

\renewcommand{\thefootnote}{$\star$}

\ShortArticleName{Monogenic Functions in Conformal Geometry}

\ArticleName{Monogenic Functions in Conformal Geometry\footnote{This paper is a
contribution to the Proceedings of the 2007 Midwest
Geometry Conference in honor of Thomas~P.\ Branson. The full collection is available at
\href{http://www.emis.de/journals/SIGMA/MGC2007.html}{http://www.emis.de/journals/SIGMA/MGC2007.html}}}

\Author{Michael EASTWOOD~$^\dag$ and John RYAN~$^\ddag$}

\AuthorNameForHeading{M.~Eastwood and J.~Ryan}

\Address{$^\dag$~Department of Mathematics, University of Adelaide,
SA 5005, Australia}
\EmailD{\href{mailto:meastwoo@member.ams.org}{meastwoo@member.ams.org}}

\Address{$^\ddag$~Department of Mathematics, University of Arkansas,
Fayetteville, AR 72701, USA}
\EmailD{\href{mailto:jryan@uark.edu}{jryan@uark.edu}}

\ArticleDates{Received August 29, 2007; Published online August 30, 2007}

\Abstract{Monogenic functions are basic to Clif\/ford analysis. On Euclidean space they are
def\/ined as smooth functions with values in the corresponding Clif\/ford algebra
satisfying a~certain system of f\/irst order dif\/ferential equations, usually referred to as
the Dirac equation. There are two equally natural extensions of these equations
to a Riemannian spin manifold only one of which is conformally invariant. We
present a straightforward exposition.}


\Classification{53A30; 58J70; 15A66}

\section{Introduction}
This article is dedicated to the memory of Tom Branson, a true gentleman,
scholar, and friend. A signif\/icant and recurring theme in his work was the
interaction between the following topics
$$\framebox{\begin{tabular}{ll}Lie Groups&Homogeneous Spaces\\
Global Symmetry&Representation Theory\end{tabular}}
\leftrightarrow
\framebox{\begin{tabular}{ll}Dif\/ferential Geometry&Local Structure\\
Local Symmetry&PDEs \end{tabular}}$$
Items in the left hand box are concerned with geometry and global analysis on a
homogeneous space $G/P$ where $G$ is a Lie group with subgroup~$P$. Items in
the right hand box are concerned with local dif\/ferential geometry `modelled' on
the homogeneous space~$G/P$. Without trying to~explain what this might mean in
general, one of the best known examples is to have conformal dif\/ferential
geometry in dimension $n\geq 3$ in the right hand box and the sphere~$S^n$,
viewed as a~homogeneous space for ${\mathrm{SO}}(n+1,1)$, in the left hand box.
The round sphere provides a~better `f\/lat model' of conformal geometry than does
Euclidean space because all the conformal Killing f\/ields integrate to genuine
conformal motions (and there is a conformal embedding
\mbox{${\mathbb{R}}^n\hookrightarrow S^n$} given by inverse stereographic projection).
Tom Branson was one of the f\/irst (e.g.~\cite{b}) systemati\-cally to exploit the
interplay between `curved' conformal dif\/ferential geometry and its f\/lat model.
He was also keen to promote the interaction between conformal geometry and the
f\/ield of Clif\/ford analysis, where {\em monogenic functions} were def\/ined and
are now studied.

Some disclaimers are in order. The results in this article are, in some sense,
already known. More precisely, there are a few mathematicians, including David
Calderbank and Vladim\'{\i}r Sou\v{c}ek, who are well aware of results in this
vein. There is no doubt that Tom Branson would have had his own distinctive
viewpoint. Our aim, therefore, is to give what we believe to be a~particularly
ef\/f\/icient formulation and method of proof with minimal prerequisites.
Motivation and consequences will be given as we go along.

Clif\/ford analysis on Euclidean space~${\mathbb{R}}^n$ is usually developed by
considering smooth functions with values in the Clif\/ford algebra
of~${\mathbb{R}}^n$. In this context one may equally well consider smooth
dif\/ferential forms on~${\mathbb{R}}^n$, i.e.\ smooth functions with values in
$\Lambda^\bullet=\bigoplus_{p=0}^n\Lambda^p$ where $\Lambda^p$ is the bundle of
$p$-forms. In fact, there is a canonical identif\/ication between
$\Lambda^\bullet{\mathbb{R}}^n$ and the Clif\/ford algebra. Under this
identif\/ication, the so-called `Dirac operator'
$D:\Lambda^\bullet\to\Lambda^\bullet$ is precisely
\[
D\omega=d\omega-d^*\omega,
\]
where $d$ is exterior derivative and $d^*$ is its adjoint. Viewed like this,
$D$ enjoys an evident extension to any Riemannian manifold where it is called
the Hodge-de~Rham operator. There is, however, an equally natural but dif\/ferent
extension of $D$ to any Riemannian manifold with spin structure. This
alternative extension turns out to be conformally invariant and hence to be
preferred as the starting point for Clif\/ford analysis on manifolds.

In this article, we shall present this alternative extension in as simple a
manner as possible. In particular, we shall not need Clif\/ford algebras for the
basic exposition, only to compare the resulting Dirac operator with the
classical one def\/ined on the usual spin bundles in~\S~\ref{classical} (and, even
then, only implicitly). The simplest exposition does not even require any
representation theory beyond the basics.

We conclude this introduction by establishing some notation and conventions. We
shall use Penrose's `abstract index notation'~\cite{OT} for tensors on a
manifold and also for representations of the orthogonal
group~${\mathrm{SO}}(n)$. Thus, we shall write $\omega_a$ to denote a smooth
$1$-form, and $X^a$ to denote a vector f\/ield. Square brackets denote skewing
over the indices they enclose~-- if \mbox{$\omega_{ab}=\omega_{[ab]}$}, then
$\omega_{ab}$ denotes a $2$-form. Similarly, we shall use round brackets to
denote symmetrisation: $\phi_{(ab)}=\frac12(\phi_{ab}+\phi_{ba})$. Repeated
indices denote the natural pairing of vectors and co-vectors:
$X\intprod\omega=X^a\omega_a$. In local co\"ordinates this is the Einstein
summation convention but the abstract index convention is itself co\"ordinate
free. Following standard practise, on a Riemannian manifold we shall write
$g_{ab}$ for the metric and $g^{ab}$ for its inverse:
$g_{ab}g^{bc}=\delta_a{}^c$ where $\delta_a{}^c$ is the Kronecker delta or
identity matrix. We shall often `raise and lower indices' without comment~-- if
$X^a$ is a vector f\/ield then $X_a\equiv g_{ab}X^b$ is the corresponding
$1$-form. Sometimes we shall write
% $\langle\alpha,\beta\rangle$ instead of $\alpha^a\beta_a$.
$\langle u,v\rangle$ instead of $u^av_a$.
We need to set some conventions for the exterior and interior product on
dif\/ferential forms~-- if $\omega_{bc\cdots d}$ has $p$ indices, then
\begin{gather}\label{extintconventions}
(v\wedge\omega)_{abc\cdots d}\equiv v_{[a}\omega_{bc\cdots d]}
\qquad\mbox{and}\qquad
(v\intprod\omega)_{c\cdots d}\equiv pv^b\omega_{bc\cdots d}.
\end{gather}
These conventions have the convenient feature that
\begin{gather}\label{convenient}
u\wedge(v\intprod\omega)+v\intprod(u\wedge\omega)=
\langle u,v\rangle\omega.\end{gather}
For further details see~\cite{w}.

\section[Conformally invariant first order operators]{Conformally invariant f\/irst order operators}
We begin with three simple examples. If $\omega_a$ is a $1$-form on a
Riemannian manifold $M$ and we rescale the metric according to
$\hat g_{ab}=\Omega^2g_{ab}$, then it is easily verif\/ied that
\[
\hat\nabla_a\omega_b=\nabla_a\omega_b-\Upsilon_a\omega_b-\Upsilon_b\omega_a
+\Upsilon^c\omega_cg_{ab}\qquad\mbox{where}\quad
\Upsilon_a=\nabla_a\Omega\big/\Omega
\]
and $\nabla_a$ is the metric connection for
$g_{ab}$ whilst $\hat\nabla_a$ is the metric connection for~$\hat g_{ab}$. It
follows that
\[
\hat\nabla_{[a}\omega_{b]}=\nabla_{[a}\omega_{b]}.
\]
This is just the conformal invariance of the exterior derivative
$d:\Lambda^1\to\Lambda^2$. More instructive examples are obtained by insisting
that when we rescale the metric so that $\hat g_{ab}=\Omega^2g_{ab}$ we also
rescale the $1$-form $\hat\omega_a=\Omega^w\omega_a$. In this case
\begin{gather*}
\hat\nabla_a\hat\omega_b=\hat\nabla_a(\Omega^w\omega_b)\\
\phantom{\hat\nabla_a\hat\omega_b}{}=(\nabla_a\Omega^w)\omega_b+\Omega^w\hat\nabla_a\omega_b\\
\phantom{\hat\nabla_a\hat\omega_b}{}= \Omega^ww\Upsilon_a\omega_b+\Omega^w\big(\nabla_a\omega_b-\Upsilon_a\omega_b
-\Upsilon_b\omega_a+\Upsilon^c\omega_cg_{ab}\big)\\
\phantom{\hat\nabla_a\hat\omega_b}{}= \Omega^w\big(\nabla_a\omega_b+(w-1)\Upsilon_a\omega_b
-\Upsilon_b\omega_a+\Upsilon^c\omega_cg_{ab}\big).
\end{gather*}
Notice that the right hand side rescales by the same power of $\Omega$ as we
decreed for our $1$-form~$\omega_a$. Thus, if we say that $\omega_a$ is
{\em conformally weighted\/} of weight $w$ instead of writing
$\hat\omega_a=\Omega^w\omega_a$, then we have shown that
\begin{gather}\label{change}\hat\nabla_a\omega_b=
\nabla_a\omega_b+(w-1)\Upsilon_a\omega_b
-\Upsilon_b\omega_a+\Upsilon^c\omega_cg_{ab}\end{gather}
for $\omega_a$ of conformal weight~$w$. More precisely, we define a conformal
manifold to be a smooth manifold equipped with an equivalence class of
Riemannian metrics under the equivalence relation $g_{ab}\mapsto\hat
g_{ab}=\Omega^2g_{ab}$ for any smooth function~$\Omega$. On such a manifold we
define a line bundle~$L$ as trivialised by any choice of metric~-- its sections
are identified as ordinary functions~$f$. For a conformally related metric
$\hat g_{ab}=\Omega^2g_{ab}$, however, the same section gives a dif\/ferent
function~$\hat f=\Omega f$. Let us write $\Lambda^0[w]$ instead of $L^w$ and
refer to its sections as {\em conformal densities\/} of weight~$w$. Similarly,
$\Lambda^p[w]=\Lambda^p\otimes L^w$ is the bundle of conformally weighted
$p$-forms of weight~$w$. For any metric $g_{ab}$ in the conformal class these
bundles acquire natural connections, namely by trivialising $L$ and employing
the metric connection. Equation (\ref{change}) now has a precise
interpretation~-- it says how the connection
$\nabla:\Lambda^1[w]\to\Lambda^1\otimes\Lambda^1[w]$ changes under a conformal
rescaling of the metric. The metric $g_{ab}$ itself acquires a tautological
interpretation as a conformally invariant section
of~$\bigodot^2\!\Lambda^1[2]$. Similarly, its inverse $g^{ab}$ has conformal
weight $-2$ and we may raise and lower indices at the expense of changing the
weight~-- if $\omega_a$ has weight~$w$, then $\omega^a$ has weight~$w-2$.

If~$w=2$, then (\ref{change}) reads
\[
\hat\nabla_a\omega_b=\nabla_a\omega_b+\Upsilon_a\omega_b-\Upsilon_b\omega_a
+\Upsilon^c\omega_cg_{ab}.
\]
It follows that
\[
\hat\nabla_{(a}\omega_{b)}-\tfrac1n\hat\nabla^c\omega_cg_{ab}=
\nabla_{(a}\omega_{b)}-\tfrac1n \nabla^c\Upsilon_cg_{ab}.
\]
Thus, we have a conformally invariant operator
\[
\textstyle\Lambda^1[2]\to\bigodot^2_\circ\!\Lambda^1[2],
\]
where $\bigodot^2_\circ\!\Lambda^1$ denote the trace-free symmetric covariant
$2$-tensors.

Taking the trace of (\ref{change}) yields
\[
\hat\nabla^a\omega_a=\nabla^a\omega_a+(n+w-2)\Upsilon^a\omega_a.
\]
We obtain, as our third example, a conformally invariant dif\/ferential operator
\[
\Lambda^1[-(n-2)]\to\Lambda^0[-n].
\]

Our fundamental theorem constructs first order conformally invariant operators
from representations of ${\mathrm{SO}}(n)$ or ${\mathrm{Spin}}(n)$. Let us
denote by~${\mathbb{W}}$, the defining representation of ${\mathrm{SO}}(n)$
on~${\mathbb{R}}^n$ and by $g_{ab}\in{\mathbb{W}}\otimes{\mathbb{W}}$ or by
$\langle\enskip,\enskip\rangle$ the standard
Euclidean metric, preserved by~${\mathrm{SO}}(n)$. We shall identify
$\Lambda^2{\mathbb{W}}$ with the Lie algebra ${\mathfrak{so}}(n)$ of
${\mathrm{SO}}(n)$ or ${\mathrm{Spin}}(n)$.
Explicitly, let us take the action
of $X_{ab}\in\Lambda^2{\mathbb{W}}$ on $w_a\in{\mathbb{W}}$ to be
given by $X_a{}^bw_b$. Also let
\[
w_a\mapsto-2g_{a[b}w_{c]}\quad\mbox{define}\quad
\iota:{\mathbb{W}}\to{\mathbb{W}}\otimes\Lambda^2{\mathbb{W}}.
\]
A representation ${\mathbb{E}}$ of ${\mathrm{SO}}(n)$ gives rise to an
irreducible tensor bundle $E$ on any oriented Riemannian $n$-manifold.
Specifically, we shall regard $E$ as induced from the orthonormal co-frame
bundle so that ${\mathbb{W}}$ gives rise to~$\Lambda^1$, the co-tangent bundle.
Similarly, if ${\mathbb{E}}$ is a ${\mathrm{Spin}}(n)$-representation, then we
find induced bundles $E$ on any oriented Riemannian spin $n$-manifold. On a
conformal manifold it is convenient, following Calderbank~\cite{c}, to induce
bundles from the orthonormal frames of~$\Lambda^1[1]$. The
reason is that $\langle\omega,\mu\rangle=g^{ab}\omega_a\mu_b$ is conformally
invariant for $\omega_a$ and $\mu_a$ of conformal weight~$1$ and so these
frames are well-defined. A further reason is that, although there is no
conformally invariant connection on~$\Lambda^1[1]$, the conformal change
(\ref{change}) is especially simple:
\begin{gather}\label{son}
\hat\nabla_a\omega_b=\nabla_a\omega_b-\Upsilon_b\omega_a
+\Upsilon^c\omega_cg_{ab}
=\nabla_a\omega_b-\Gamma_{ab}{}^c\omega_c
\qquad\mbox{where}\quad
\Gamma_{abc}=-2g_{a[b}\Upsilon_{c]}.\end{gather}
We shall denote by $E$ the bundle induced on an oriented conformal manifold from a
representation ${\mathbb{E}}$ of ${\mathrm{SO}}(n)$. {From} the defining
representation ${\mathbb{W}}$ of ${\mathrm{SO}}(n)$ we obtain the bundle
$\Lambda^1[1]$.

\begin{theorem}\label{fundamental} Suppose that ${\mathbb{E}}$ and
${\mathbb{F}}$ are representations of\/ ${\mathrm{SO}}(n)$ and that
\[
\pi:{\mathbb{W}}\otimes{\mathbb{E}}\to{\mathbb{F}}
\]
is a homomorphism of ${\mathrm{SO}}(n)$-modules. Let
$\rho:\Lambda^2{\mathbb{W}}\otimes{\mathbb{E}}\to{\mathbb{E}}$ denote the
action of\/ ${\mathfrak{so}}(n)$ on~${\mathbb{E}}$. Suppose that the
composition
\[
{\mathbb{W}}\otimes{\mathbb{E}}\xrightarrow{\iota\otimes{\mathrm{Id}}}
{\mathbb{W}}\otimes\Lambda^2{\mathbb{W}}\otimes{\mathbb{E}}
\xrightarrow{{\mathrm{Id}}\otimes\rho}{\mathbb{W}}\otimes{\mathbb{E}}
\xrightarrow{\pi}{\mathbb{F}}
\]
is equal to $w\pi:{\mathbb{W}}\otimes{\mathbb{E}}\to{\mathbb{F}}$ for some
constant~$w$. Then there is a conformally invariant first order linear
differential operator
\[
D:E[w]\to F[w-1]
\]
whose symbol is induced by~$\pi$. If, instead, ${\mathbb{E}}$ and
${\mathbb{F}}$ are representations of\/ ${\mathrm{Spin}}(n)$, then we obtain
the corresponding conclusion on any Riemannian spin manifold.
\end{theorem}

\begin{proof}
Recall the formula (\ref{son}) for a change of connection on $\Lambda^1[1]$.
Therefore, for sections $\phi$ of the associated bundle $E$,
\[
\hat\nabla\phi=\nabla\phi-\Gamma\phi
\]
where $\Gamma\phi$ is the image of $\Upsilon\otimes\phi$ under the composition
\[
\Lambda^1\otimes E\xrightarrow{\iota\otimes{\mathrm{Id}}}
\Lambda^1\otimes\Lambda^2\otimes E\xrightarrow{{\mathrm{Id}}\otimes\rho}
\Lambda^1\otimes E.
\]
For sections of $E[w]$, therefore,
\[
\hat\nabla\phi=\nabla\phi+w\Upsilon\otimes\phi-\Gamma\phi.
\]
Hence, for $D=\pi\nabla$ we have
\[
\hat D\phi=D\phi+w\pi(\Upsilon\otimes\phi)-\pi\Gamma\phi=D\phi,
\]
as required.
\end{proof}

The three conformally invariant operators given at the beginning of this section
are easily obtained from Theorem~\ref{fundamental}. Specifically, the
irreducible decomposition
\[
\textstyle{\mathbb{W}}\otimes{\mathbb{W}}=\Lambda^2{\mathbb{W}}\oplus
\bigodot^2_\circ\!{\mathbb{W}}\oplus{\mathbb{R}}
\]
is given explicitly by
\[
\textstyle\xi_a\omega_b=\xi_{[a}\omega_{b]}+
\big(\xi_{(a}\omega_{b)}-\frac1n\xi^c\omega_cg_{ab}\big)
+\frac1n\xi^c\omega_cg_{ab}
\]
from which we obtain three possible symbols with ${\mathbb{E}}={\mathbb{W}}$ so
that $E=\Lambda^1[1]$. The composition
\[
{\mathbb{W}}\otimes{\mathbb{E}}\xrightarrow{\iota\otimes{\mathrm{Id}}}
{\mathbb{W}}\otimes\Lambda^2{\mathbb{W}}\otimes{\mathbb{E}}
\xrightarrow{{\mathrm{Id}}\otimes\rho}{\mathbb{W}}\otimes{\mathbb{E}}
\]
is
$\xi_a\omega_b\mapsto-2g_{a[b}\xi_{c]}\omega_d\mapsto
\xi_b\omega_a-g_{ab}\xi^c\omega_c\equiv\sigma_{ab}$ and to apply
Theorem~\ref{fundamental} we must see whether $\pi\sigma$ is a multiple of
$\pi(\xi\otimes\omega)$. We find
\begin{alignat*}{4}
& \xi_a\omega_b\stackrel{\pi}{\longmapsto}\xi_{[a}\omega_{b]}&&\mbox{gives}&&
\sigma_{[ab]}=-\xi_{[a}\omega_{b]},&\\
&\xi_a\omega_b\stackrel{\pi}{\longmapsto}
\xi_{(a}\omega_{b)}-\tfrac1n\xi^c\omega_cg_{ab}\quad &&\mbox{gives}\quad &&
\sigma_{(ab)}-\tfrac1n\sigma^c{}_cg_{ab}
=\big(\xi_{(a}\omega_{b)}-\tfrac1n\xi^c\omega_cg_{ab}\big),&\\
&\xi_a\omega_b\stackrel{\pi}{\longmapsto}\xi^a\omega_a&&\mbox{gives}&&
\sigma^a{}_a=-(n-1)\xi^a\omega_a.&
\end{alignat*}
We obtain our three basic conformally invariant operators
\[
\begin{tabular}{c|c|c|c|c}
${\mathbb{F}}$&$w$&$F$&$E[w]$&$F[w-1]$\\ \hline\hline
\vphantom{\Big(}$\Lambda^2{\mathbb{W}}$&$-1$&$\Lambda^2[2]$&$\Lambda^1$&
$\Lambda^2$\\ \hline
\vphantom{\Big(}$\bigodot_\circ^2\!{\mathbb{W}}$&$1$&
$\bigodot_\circ^2\!\Lambda^1[2]$&$\Lambda^1[2]$&
$\bigodot_\circ^2\!\Lambda^1[2]$\\ \hline
\vphantom{\Big(}${\mathbb{R}}$&$-(n-1)$&$\Lambda^0$&$\Lambda^1[-(n-2)]$&
$\Lambda^0[-n]$
\end{tabular}
\]
in accordance with Theorem~\ref{fundamental}.

In stating Theorem~\ref{fundamental} we were imprecise concerning whether
${\mathbb{E}}$ and ${\mathbb{F}}$ should be real or complex representations of
${\mathrm{SO}}(n)$ or ${\mathrm{Spin}}(n)$. In fact, there are two
versions~-- if the representations are real (as in our three examples), then the
corresponding induced bundles are real but if the representations are complex
(as is more usual in representation theory) then the bundles are complex.

The hypotheses of Theorem~\ref{fundamental} are automatically satisfied if
${\mathbb{E}}$ is an irreducible complex representation of ${\mathrm{SO}}(n)$
or ${\mathrm{Spin}}(n)$. In this case the decomposition of
${\mathbb{W}}\otimes{\mathbb{E}}$ into irreducibles is multiplicity-free. Thus,
if $\pi:{\mathbb{W}}\otimes{\mathbb{E}}\to{\mathbb{F}}$ is projection onto any
of the irreducibles, then any other ${\mathfrak{so}}(n)$-invariant homomorphism
${\mathbb{W}}\otimes{\mathbb{E}}\to{\mathbb{E}}$ must be a multiple thereof. In
this case, it is just a matter of identifying the constant $w$ in order to
classify the first order conformally invariant dif\/ferential operators between
irreducible tensor or spinor bundles. This is precisely what Fegan does
in~\cite{hf}, using Casimir operators to compute~$w$. In~\cite{css3}, \v{C}ap,
Slov\'ak, and Sou\v{c}ek extended Fegan's method to higher order
operators.

\section{Monogenic functions and Dirac operators}\label{defofD}
The construction and conformal invariance of monogenic functions also comes
from Theorem~\ref{fundamental}. In distinction to the discussion at the end of
the previous section, the representation ${\mathbb{E}}$ will be real and
reducible. As a vector space, let ${\mathbb{E}}=\Lambda^\bullet{\mathbb{W}}$
where, as usual, ${\mathbb{W}}$ is the (real) defining representation of
${\mathrm{SO}}(n)$. Define a homomorphism
\begin{gather}\label{defdot}
{\mathbb{W}}\otimes{\mathbb{E}}\xrightarrow{\epsilon}{\mathbb{E}}\qquad
\mbox{by}\quad v\otimes e\mapsto v.e\equiv
v\wedge e-v\intprod e.
\end{gather}
There are two natural ways in which ${\mathbb{E}}$ is an
${\mathfrak{so}}(n)$-module. The more obvious corresponds to the
${\mathrm{SO}}(n)$-action defined by that on~${\mathbb{W}}$. Evidently,
$\epsilon$ is a homomorphism of ${\mathrm{SO}}(n)$-modules. The geometric
consequence of this is a dif\/ferential operator from $E$ to itself defined on
any Riemannian manifold as the composition
\[
E\xrightarrow{\nabla}\Lambda^1\otimes E\xrightarrow{\epsilon}E,
\]
where $\nabla$ is induced by the metric connection on the co-frame bundle. This
is the Hodge--de~Rham operator. For this choice of
$\rho:\Lambda^2{\mathbb{W}}\otimes{\mathbb{E}}\to{\mathbb{E}}$, however,
\[
{\mathbb{W}}\otimes{\mathbb{E}}\xrightarrow{\iota\otimes{\mathrm{Id}}}
{\mathbb{W}}\otimes\Lambda^2{\mathbb{W}}\otimes{\mathbb{E}}
\xrightarrow{{\mathrm{Id}}\otimes\rho}{\mathbb{W}}\otimes{\mathbb{E}}
\xrightarrow{\epsilon}{\mathbb{E}}
\]
is not a multiple of~$\epsilon:{\mathbb{W}}\otimes{\mathbb{E}}\to{\mathbb{E}}$.
The Hodge--de~Rham operator is not conformally invariant.

An alternative but equally natural way in which ${\mathbb{E}}$ is an
${\mathfrak{so}}(n)$-module is given by
\begin{gather}\label{spinaction}\textstyle
{\mathfrak{so}}(n)\otimes{\mathbb{E}}=\Lambda^2{\mathbb{W}}\otimes{\mathbb{E}}
\ni v\wedge w\otimes e\stackrel{\sigma}{\longmapsto}
-\frac18(v.w.e-w.v.e)\in{\mathbb{E}}.
\end{gather}
The constant $-1/8$ ensures that this is, indeed, a representation
of~${\mathfrak{so}}(n)$. To verify this and other properties of this
construction, the following lemmata are useful.

\begin{lemma}\label{cliffordmodule}
% For all $v,w\in{\mathbb{W}}$ and $e\in{\mathbb{E}}$ we have
% \begin{gather}\label{cliff}v.w.e+w.v.e=-2\langle v,w\rangle e.\end{gather}
For all $u,v\in{\mathbb{W}}$ and $e\in{\mathbb{E}}$ we have
\begin{gather}\label{cliff}u.v.e+v.u.e=-2\langle u,v\rangle e.\end{gather}
\end{lemma}

\begin{proof} A simple computation from (\ref{convenient}).
\end{proof}
\begin{lemma} For all $u,v,w\in{\mathbb{W}}$ and $e\in{\mathbb{E}}$ we have
\begin{gather}\label{key}u.v.w.e-v.w.u.e-u.w.v.e+w.v.u.e=
-4\langle u,v\rangle w.e+4\langle u,w\rangle v.e.\end{gather}
\end{lemma}
\begin{proof} From (\ref{cliff}) we find
\begin{gather*}
u.v.w.e-v.w.u.e =u.v.w.e+v.u.w.e-v.u.w.e-v.w.u.e\\
\phantom{u.v.w.e-v.w.u.e}{} =-2\langle u,v\rangle w.e+2\langle u,w\rangle v.e.
\end{gather*}
The other terms in (\ref{key}) are dealt with similarly.
\end{proof}
\begin{proposition}\label{first}
The action \eqref{spinaction} makes ${\mathbb{E}}$ into a representation
of\/ ${\mathfrak{so}}(n)$. With this structure
% $\sigma:\Lambda^2{\mathbb{W}}\otimes{\mathbb{E}}\to{\mathbb{E}}$
$\epsilon:{\mathbb{W}}\otimes{\mathbb{E}}\to{\mathbb{E}}$
is a homomorphism of\/ ${\mathfrak{so}}(n)$-modules.
\end{proposition}

\begin{proof} On simple vectors, the Lie bracket on
${\mathfrak{so}}(n)=\Lambda^2{\mathbb{W}}$ is given by
\[
\textstyle [t\wedge u,v\wedge w]=
\frac12\big(\langle u,v\rangle t\wedge w-\langle u,w\rangle t\wedge v-
\langle t,v\rangle u\wedge w+\langle t,w\rangle u\wedge v\big).
\]
The assertions are straightforward calculations on simple vectors
using~(\ref{key}).
\end{proof}

\begin{lemma}\label{metricimage}
If we define
${\mathbb{W}}\otimes{\mathbb{W}}\otimes{\mathbb{E}}\to{\mathbb{E}}$ by
$v\otimes w\otimes e\mapsto v.w.e$, then $g_{ab}\otimes e\mapsto -ne$.
\end{lemma}
\begin{proof} If $u_1,u_2,\ldots,u_n$ is an orthonormal basis
of~${\mathbb{W}}$, then
$g_{ab}=\sum\limits_{i,j=1}^nu_i\otimes u_j$. The result follows immediately
from~(\ref{cliff}).
\end{proof}

\begin{proposition}\label{second} The composition
\begin{gather}\label{whatwewant}
{\mathbb{W}}\otimes{\mathbb{E}}\xrightarrow{\iota\otimes{\mathrm{Id}}}
{\mathbb{W}}\otimes\Lambda^2{\mathbb{W}}\otimes{\mathbb{E}}
\xrightarrow{{\mathrm{Id}}\otimes\sigma}{\mathbb{W}}\otimes{\mathbb{E}}
\xrightarrow{\epsilon}{\mathbb{E}}\end{gather}
is equal to
$-\frac{n-1}2\epsilon:{\mathbb{W}}\otimes{\mathbb{E}}\to{\mathbb{E}}$.
\end{proposition}
\begin{proof} Let us perform this computation using abstract indices. For this
purpose, if we write $\epsilon:{\mathbb{W}}\otimes{\mathbb{E}}\to{\mathbb{E}}$
as $w_a\otimes e_\alpha\mapsto w_a\epsilon^a{}_\alpha{}^\beta e_\beta$, then
Lemma~\ref{metricimage} says that
\begin{gather}\label{MIwithindices}
\epsilon^a{}_\alpha{}^\beta\epsilon_{a\beta}{}^\gamma=
-n\delta_\alpha{}^\gamma.\end{gather}
Let us write $\sigma:\Lambda^2{\mathbb{W}}\otimes{\mathbb{E}}\to{\mathbb{E}}$
as $X_{ab}\otimes e_\alpha\mapsto X_{ab}\sigma^{ab}{}_\alpha{}^\beta e_\beta$,
where
\[
\textstyle\sigma^{ab}{}_\alpha{}^\gamma=
-\frac18(\epsilon^a{}_\alpha{}^\beta\epsilon^b{}_\beta{}^\gamma-
\epsilon^b{}_\alpha{}^\beta\epsilon^a{}_\beta{}^\gamma).
\]
Then (\ref{whatwewant}) becomes
\[
w_ae_\alpha\mapsto\epsilon^a{}_\alpha{}^\beta
(-2g_{a[b}w_{c]}\sigma^{bc}{}_\beta{}^\gamma e_\gamma)=
-2\epsilon^a{}_\alpha{}^\beta
w_c\sigma_a{}^c{}_\beta{}^\gamma e_\gamma
\]
and so we are required to show that
\[
\textstyle -2\epsilon^a{}_\alpha{}^\beta
\sigma_a{}^c{}_\beta{}^\gamma=
-\frac{n-1}2\epsilon^c{}_\alpha{}^\gamma\quad\mbox{i.e.}\quad
4\epsilon^a{}_\alpha{}^\beta\sigma_a{}^c{}_\beta{}^\gamma=
(n-1)\epsilon^c{}_\alpha{}^\gamma.
\]
Using (\ref{MIwithindices}), we compute:
\[
\textstyle 4\epsilon^a{}_\alpha{}^\beta\sigma_a{}^c{}_\beta{}^\gamma=
-\frac12\epsilon^a{}_\alpha{}^\beta
(\epsilon_{a\beta}{}^\eta\epsilon^c{}_\eta{}^\gamma-
\epsilon^c{}_\beta{}^\eta\epsilon_{a\eta}{}^\gamma)=
\frac{n}2\epsilon^c{}_\alpha{}^\gamma
+\frac12\epsilon^a{}_\alpha{}^\beta
\epsilon^c{}_\beta{}^\eta\epsilon_{a\eta}{}^\gamma.
\]
To continue, let us write (\ref{key}) with abstract indices:
\[
\epsilon^a{}_\alpha{}^\beta\epsilon^b{}_\beta{}^\eta\epsilon^c{}_\eta{}^\gamma
-\epsilon^b{}_\alpha{}^\beta\epsilon^c{}_\beta{}^\eta\epsilon^a{}_\eta{}^\gamma
-\epsilon^a{}_\alpha{}^\beta\epsilon^c{}_\beta{}^\eta\epsilon^b{}_\eta{}^\gamma
+\epsilon^c{}_\alpha{}^\beta\epsilon^b{}_\beta{}^\eta\epsilon^a{}_\eta{}^\gamma
=-4g^{ab}\epsilon^c{}_\alpha{}^\gamma+4g^{ac}\epsilon^b{}_\alpha{}^\gamma.
\]
Tracing over $a$ and $b$ and using (\ref{MIwithindices}) gives
\[
\textstyle
-n\epsilon^c{}_\alpha{}^\gamma
-\epsilon^a{}_\alpha{}^\beta\epsilon^c{}_\beta{}^\eta\epsilon_{a\eta}{}^\gamma
-\epsilon^a{}_\alpha{}^\beta\epsilon^c{}_\beta{}^\eta\epsilon_{a\eta}{}^\gamma
-n\epsilon^a{}_\alpha{}^\gamma
=-4n\epsilon^c{}_\alpha{}^\gamma+4\epsilon^c{}_\alpha{}^\gamma.
\]
If follows that
\[
\epsilon^a{}_\alpha{}^\beta\epsilon^c{}_\beta{}^\eta\epsilon_{a\eta}{}^\gamma
=(n-2)\epsilon^c{}_\alpha{}^\gamma
\]
and hence that
\[
\textstyle 4\epsilon^a{}_\alpha{}^\beta\sigma_a{}^c{}_\beta{}^\gamma=
\frac{n}2\epsilon^c{}_\alpha{}^\gamma
+\frac{n-2}2\epsilon^c{}_\alpha{}^\gamma=
(n-1)\epsilon^c{}_\alpha{}^\gamma,
\]
as required.
\end{proof}

\begin{theorem}\label{monogenicinvariance}
Let ${\mathbb{E}}$ denote that representation of\/ ${\mathrm{Spin}}(n)$
corresponding to~\eqref{spinaction}. Let $E$ denote the corresponding
bundle on a conformal spin manifold induced from the orthonormal frames
of~$\Lambda^1[1]$. Then there is a conformally invariant first order linear
differential operator
\[
\textstyle E[-\frac{n-1}2]\xrightarrow{D}E[-\frac{n+1}2]
\]
whose symbol is induced by
$\epsilon:{\mathbb{W}}\otimes{\mathbb{E}}\to{\mathbb{E}}$.
\end{theorem}

\begin{proof}
Proposition~\ref{first} ensures that the statement of this theorem makes sense
and now Proposition~\ref{second} ensures that the criterion of
Theorem~\ref{fundamental} is satisfied.
\end{proof}

\begin{remark} This operator $D$ extends the Dirac operator from Clif\/ford
analysis on ${\mathbb{R}}^n$ to a~general spin manifold. Elements of its kernel
are referred to as `monogenic functions'.
\end{remark}

\begin{remark} Almost no representation theory is needed here~-- we have only
used that ${\mathrm{Spin}}(n)$ is the simply-connected connected Lie group
whose Lie algebra is ${\mathfrak{so}}(n)$. In fact, we shall see in the next
section that spin is essential~-- the representation (\ref{spinaction}) of
${\mathfrak{so}}(n)$ does not arise from an action of~${\mathrm{SO}}(n)$. The
proof of Theorem~\ref{monogenicinvariance} is purely computational and the only
ingredients in this computation are (\ref{key}) and (\ref{MIwithindices}) both
of which follow easily from Lemma~\ref{cliffordmodule}. In fact, the particular
constant in Lemma~\ref{metricimage} (leading to~(\ref{MIwithindices})) is
unimportant. We would obtain the same result (with the same conformal weight)
from $g_{ab}\otimes e\mapsto\kappa e$ for any non-zero~$\kappa$.
\end{remark}

\section[Clifford algebras and the classical Dirac operator]{Clif\/ford algebras and the classical Dirac operator}\label{classical}

The Dirac operator from Clif\/ford analysis is usually introduced via the
Clif\/ford algebra~${\mathcal{C}}\ell({\mathbb{W}})$, defined as the tensor
algebra $\bigotimes^\bullet{\mathbb{W}}$ modulo the two-sided ideal generated
by
\[
v\otimes w+w\otimes v+2\langle v,w\rangle\qquad\forall\,v,w\in{\mathbb{W}}.
\]
Writing the multiplication in ${\mathcal{C}}\ell({\mathbb{W}})$ as
juxtaposition, there is a canonical identification
\begin{gather}\label{cliff-ext}\textstyle
{\mathbb{E}}\equiv\Lambda^\bullet{\mathbb{W}}\xrightarrow{\simeq\enskip}
{\mathcal{C}}\ell({\mathbb{W}})\qquad\mbox{by}\quad v\wedge w\mapsto
\frac12(vw-wv).
\end{gather}
The conclusion (\ref{cliff}) of Lemma~\ref{cliffordmodule} is precisely that
$\epsilon:{\mathbb{W}}\otimes{\mathbb{E}}\to{\mathbb{E}}$ extends to a
representation ${\mathcal{C}}\ell({\mathbb{W}})$, namely to a homomorphism of
algebras ${\mathcal{C}}\ell({\mathbb{W}})\to{\mathrm{End}}({\mathbb{E}})$.
Equivalently, if we transport the Clif\/ford algebra structure from
${\mathcal{C}}\ell({\mathbb{W}})$ to $\Lambda^\bullet{\mathbb{W}}$ using
(\ref{cliff-ext}) then (\ref{defdot}) gives the Clif\/ford product between
${\mathbb{W}}$ and~$\Lambda^\bullet{\mathbb{W}}$ and so the action of
${\mathcal{C}}\ell({\mathbb{W}})$ on ${\mathbb{E}}$ induced by $\epsilon$
becomes the action of ${\mathcal{C}}\ell({\mathbb{W}})$ on itself by left
multiplication. For further details see~\cite{h}.

A standard rationale for introducing the Clif\/ford algebra is in providing a
concrete realisation of~${\mathrm{Spin}}(n)$, namely as a subgroup of the group
of invertible elements in~${\mathcal{C}}\ell({\mathbb{R}}^n)$. Then it is clear
that ${\mathbb{E}}$ under (\ref{cliff-ext}) is a real representation
of~${\mathrm{Spin}}(n)$~-- simply restrict the action of
${\mathcal{C}}\ell({\mathbb{R}}^n)$ to~${\mathrm{Spin}}(n)$.

The basic spin representations of ${\mathfrak{so}}(n)$ are complex
representations. Hence, in order to make contact with ${\mathbb{E}}$ as an
${\mathfrak{so}}(n)$-module given by (\ref{spinaction}), it is necessary to
complexify. Let us suppose that $n$ is even. Then we may write ${\mathbb{CW}}$
as the direct sum of two totally null subspaces, using the complexified metric
as a dual pairing between them. Specifically, let us take
\begin{gather}\label{innerproduct}
{\mathbb{CW}}={\mathbb{U}}\oplus{\mathbb{U}}^*\ni\alpha+\beta
\qquad\mbox{so that}\quad \|\alpha+\beta\|^2=2\alpha\intprod\beta.\end{gather}
Let use abstract indices to write elements of ${\mathbb{U}}$ as $\alpha_a$ and
elements of ${\mathbb{U}}^*$ as~$\beta^a$, without being alarmed that the index
$a$ now runs only over half the range that it did in previous sections. Forms
now decompose according to `type'. Specifically,
${\mathbb{CE}}=\Lambda^\bullet{\mathbb{CW}}$ and
\[
\Lambda^r{\mathbb{CW}}=
\bigoplus_{p+q=r}\Lambda^p{\mathbb{U}}\otimes\Lambda^q{\mathbb{U}}^*.
\]
Splitting the formulae (\ref{extintconventions}) according to this
decomposition, we find
\begin{gather*}
\mbox{for }\Lambda^p{\mathbb{U}}\otimes\Lambda^q{\mathbb{U}}^*\ni
\omega\leftrightsquigarrow
\omega_{\mbox{\scriptsize$\underbrace{bc\cdots d}_p$}}
{}^{\mbox{\scriptsize$\overbrace{fg\cdots h}^q$}},\\
\qquad \left\lbrace\!\!
\begin{array}l
(\alpha+\beta)\wedge\omega\leftrightsquigarrow
\alpha_{[a}\omega_{bc\cdots d]}{}^{fg\cdots h}+
(-1)^p\beta^{[e}\omega_{bc\cdots d}{}^{fg\cdots h]},\vspace{2mm}\\
(\alpha+\beta)\intprod\omega\leftrightsquigarrow
p\beta^b\omega_{bc\cdots d}{}^{fg\cdots h}
+(-1)^pq\alpha_f\omega_{bc\cdots d}{}^{fg\cdots h}.
\end{array}\right.
\end{gather*}
Consequently, as a complexification of the action of ${\mathbb{W}}$ on
${\mathbb{E}}=\Lambda^\bullet{\mathbb{W}}$ given by (\ref{defdot}) we obtain an
action of ${\mathbb{CW}}$ on
$\Lambda^\bullet{\mathbb{U}}\otimes\Lambda^\bullet{\mathbb{U}}^*$ given by
\begin{gather*}
(\alpha+\beta).\omega\leftrightsquigarrow
\alpha_{[a}\omega_{bc\cdots d]}{}^{fg\cdots h}+
(-1)^p\beta^{[e}\omega_{bc\cdots d}{}^{fg\cdots h]}
-p\beta^b\omega_{bc\cdots d}{}^{fg\cdots h}
-(-1)^pq\alpha_f\omega_{bc\cdots d}{}^{fg\cdots h}.
\end{gather*}
In accordance with (\ref{cliff}) and~(\ref{innerproduct}), or as may be
verified by direct computation, this action has the property that
\[
(\alpha+\beta).(\alpha+\beta).\omega=-2\alpha_a\beta^a\omega.
\]
There is, however, another action with this property. Specifically, it is
easily verified that
\[
(\alpha+\beta){:}\omega\leftrightsquigarrow
\sqrt{2}\left(\alpha_{[a}\omega_{bc\cdots d]}{}^{fg\cdots h}
-p\beta^b\omega_{bc\cdots d}{}^{fg\cdots h}\right)\implies
(\alpha+\beta){:}(\alpha+\beta){:}\omega=-2\alpha_a\beta^a\omega.
\]

\begin{proposition}\label{Phi}
There is an automorphism of
$\Lambda^\bullet{\mathbb{U}}\otimes\Lambda^\bullet{\mathbb{U}}^*$ that converts
$(\alpha+\beta).\omega$ into $(\alpha+\beta){:}\omega$ for all
$\alpha+\beta\in{\mathbb{CW}}$.
\end{proposition}

\begin{proof}More precisely, we want to find
$\Phi:\Lambda^\bullet{\mathbb{U}}\otimes\Lambda^\bullet{\mathbb{U}}^*\to
\Lambda^\bullet{\mathbb{U}}\otimes\Lambda^\bullet{\mathbb{U}}^*$, an
invertible linear transformation, so that
\begin{gather}\label{desiredproperty}
(\alpha+\beta){:}\Phi(\omega)=\Phi((\alpha+\beta).\omega)\qquad\!\!
\forall\:(\alpha+\beta)\in{\mathbb{U}}\oplus{\mathbb{U}}^*={\mathbb{CW}}
\quad \mbox{and}\quad
\omega\in\Lambda^\bullet{\mathbb{U}}\otimes\Lambda^\bullet{\mathbb{U}}^*.\!\!
\end{gather}
Let us take
\[
\Phi(1)\equiv 1+\delta_b{}^f+\delta_{[b}{}^{[f}\delta_{c]}{}^{g]}+\cdots
+\delta_{[b}{}^{[f}\delta_c{}^g\cdots\delta_{d]}{}^{h]}+\cdots.
\]
Then, in order for (\ref{desiredproperty}) to hold, we should have
$\alpha{:}\Phi(1)=\Phi(\alpha.1)=\Phi(\alpha)$ whence
\[
\Phi(\alpha_b)=\sqrt{2}\left(\alpha_b+\alpha_{[b}\delta_{c]}{}^f
+\alpha_{[b}\delta_c{}^f\delta_{d]}{}^g+\cdots\right).
\]
Similarly, we should have $\beta{:}\Phi(1)=\Phi(\beta.1)=\Phi(\beta)$, which
forces
\[
\Phi(\beta^f)=-\sqrt{2}\left(\beta^f+2\beta^{[f}\delta_b{}^{g]}
+3\beta^{[f}\delta_b{}^g\delta_c{}^{h]}+\cdots\right).
\]
Now, however, there is something to check because
$\alpha.\beta+\beta.\alpha=-2\alpha_a\beta^a$ and so for consistency it must
be that
\begin{gather}\label{check}
\alpha{:}\Phi(\beta)+\beta{:}\Phi(\alpha)=-2\alpha_a\beta^a\Phi(1).
\end{gather}
This is readily verified as follows
\begin{gather*}
\alpha_b{:}\Phi(\beta^f)=-2\left(\alpha_b\beta^f
+2\alpha_{[b}\beta^{[f}\delta_{c]}{}^{g]}
+\cdots\right),\\
\beta^f{:}\Phi(\alpha_b)=-2\left(\beta^a\alpha_a
+2\beta^a\alpha_{[a}\delta_{b]}{}^f
+3\beta^a\alpha_{[a}\delta_b{}^f\delta_{c]}{}^g+\cdots\right)\\
\phantom{\beta^f{:}\Phi(\alpha_b)}{}=-2\left(\beta^a\alpha_a
+\beta^a\alpha_a\delta_b{}^f-\beta^a\alpha_b\delta_a{}^f
+\beta^a\alpha_a\delta_{[b}{}^f\delta_{c]}{}^g
-2\alpha_{[b}\beta{}^{[f}\delta_{c]}{}^{g]}+\cdots\right)\\
\implies\alpha_b{:}\Phi(\beta^f)+\beta^f{:}\Phi(\alpha_b)=
-2\beta^a\alpha_a
\left(1+\delta_b{}^f+\delta_{[b}{}^f\delta_{c]}{}^g+\cdots\right)=
-2\alpha_a\beta^a\Phi(1).
\end{gather*}
So far, we know $\Phi$ on ${\mathbb{C}}\oplus{\mathbb{CW}}$. But now we may use
the desired property
\[
\Phi(v.w)=v{:}\Phi(w)\qquad\mbox{for} \ \ v,w\in{\mathbb{CW}}
\]
of $\Phi$ to extend its definition to
${\mathbb{C}}\oplus{\mathbb{CW}}\oplus\Lambda^2{\mathbb{CW}}$: it is only
necessary for consistency to check~(\ref{check}) and, similarly, that
\[
\alpha{:}\Phi(\alpha)=0\qquad\mbox{and}\qquad\beta{:}\Phi(\beta)=0.
\]
Next one uses the desired property (\ref{desiredproperty}) to extend the
definition of $\Phi$ to $\bigoplus_{r=0}^3\Lambda^r{\mathbb{CW}}$ and so on by
induction. The details are left to the reader.
\end{proof}
\begin{remark}For those readers who know the theory of Clif\/ford algebras, we
have in $\Lambda^\bullet{\mathbb{CW}}$ and
$\Lambda^\bullet{\mathbb{U}}\otimes\Lambda^\bullet{\mathbb{U}}^*=
{\mathrm{End}}({\mathbb{U}})$ two dif\/ferent realisations of the Clif\/ford
algebra~${\mathcal{C}}\ell({\mathbb{CW}})$. The mapping~$\Phi$ constructed in
the proof of Proposition~\ref{Phi} is the unique unital isomorphism between
them.
\end{remark}
The geometric consequences of Proposition~\ref{Phi} are as follows. Recall that
the Dirac operator $D$ in the sense of Clif\/ford analysis was built in
\S~\ref{defofD} from (\ref{defdot}) and (\ref{spinaction}) both of which are
defined in terms of $v.e$ for $v\in{\mathbb{W}}$ and
$e\in{\mathbb{E}}\equiv\Lambda^\bullet{\mathbb{W}}$. We have just seen that, in
case $n$ is even, we can write
\[
{\mathbb{CE}}=\Lambda^\bullet{\mathbb{CW}}
=\Lambda^\bullet({\mathbb{U}}\oplus{\mathbb{U}}^*)=
\Lambda^\bullet{\mathbb{U}}\otimes\Lambda^\bullet{\mathbb{U}}^*.
\]
Therefore, Proposition~\ref{Phi} may be viewed as providing an automorphism
$\Phi$ of ${\mathbb{CE}}$ so that
\[
\Phi(v.e)=v{:}\Phi(e)\qquad\forall\; v\in{\mathbb{CW}}
\quad \mbox{and} \quad e\in{\mathbb{CE}}.
\]
In particular, this is true for $v\in{\mathbb{W}}\hookrightarrow{\mathbb{CW}}$.
Though it is essential to complexify ${\mathbb{W}}$ and also to choose a
splitting ${\mathbb{CW}}={\mathbb{U}}\oplus{\mathbb{U}}^*$ in order to define
$\Phi$ as a complex linear automorphism of ${\mathbb{CE}}$ and though it is
also necessary in order to write down the formula
\begin{gather}\label{defcolon}
(\alpha+\beta){:}\omega\leftrightsquigarrow
\sqrt{2}\left(\alpha_{[a}\omega_{bc\cdots d]}{}^{fg\cdots h}
-p\beta^b\omega_{bc\cdots d}{}^{fg\cdots h}\right)\qquad\!\!\mbox{for}\ \
\left\lbrace\!\!\begin{array}l
\alpha+\beta\in{\mathbb{U}}\oplus{\mathbb{U}}^*={\mathbb{CW}},\\
\omega\in\Lambda^\bullet{\mathbb{U}}\otimes\Lambda^\bullet{\mathbb{U}}^*
={\mathbb{CE}},\end{array}\right.\!\!\!\!
\end{gather}
it is not necessary for $v$ to be complex in order that $v{:}\omega$ be
perfectly well-defined. The key point to observe about~(\ref{defcolon}) is
that, when viewed on
$\Lambda^\bullet{\mathbb{U}}\otimes\Lambda^\bullet{\mathbb{U}}^*$, the action
$\omega\mapsto v{:}\omega$ is entirely on~$\Lambda^\bullet{\mathbb{U}}$ with
$\Lambda^\bullet{\mathbb{U}}^*$ as a passenger. It follows that the
complexification of the representation (\ref{spinaction}) of~${\mathfrak{so}}(n)$ is isomorphic to
${\mathbb{S}}\otimes{\Lambda^\bullet\mathbb{U}}^*$ where ${\mathbb{S}}$ is
$\Lambda^\bullet{\mathbb{U}}$ regarded as an ${\mathfrak{so}}(n)$-module
according to
\[
\textstyle{\mathfrak{so}}(n)\otimes{\mathbb{S}}=
\Lambda^2{\mathbb{W}}\otimes\Lambda^\bullet{\mathbb{U}} \ni v\wedge
w\otimes\omega\mapsto-\frac18(v{:}w{:}\omega-w{:}v{:}\omega)\in
{\mathbb{S}}.
\]
This action manifestly preserves the splitting
$\Lambda^\bullet{\mathbb{U}}=
\Lambda^{\mathrm{even}}{\mathbb{U}}\oplus\Lambda^{\mathrm{odd}}{\mathbb{U}}$.
It is the standard spin representation
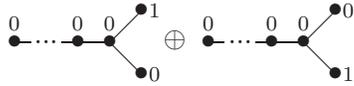
\[\hspace*{-1mm}
\raisebox{-10pt}{\begin{picture}(60,24)(6,-7)
\put(12,5){\line(1,0){6}}
\put(30,5){\line(1,0){18}}
\put(49,6){\line(1,1){10}}
\put(49,4){\line(1,-1){10}}
\put(21,5){\makebox(0,0){$.$}}
\put(24,5){\makebox(0,0){$.$}}
\put(27,5){\makebox(0,0){$.$}}
\put(12,5){\makebox(0,0){$\bullet$}}
\put(36,5){\makebox(0,0){$\bullet$}}
\put(48,5){\makebox(0,0){$\bullet$}}
\put(60,17){\makebox(0,0){$\bullet$}}
\put(60,-7){\makebox(0,0){$\bullet$}}
\put(12,12){\makebox(0,0){\scriptsize$0$}}
\put(36,12){\makebox(0,0){\scriptsize$0$}}
\put(48,12){\makebox(0,0){\scriptsize$0$}}
\put(65,17){\makebox(0,0){\scriptsize$1$}}
\put(65,-7){\makebox(0,0){\scriptsize$0$}}
\end{picture}}
\oplus
\raisebox{-10pt}{\begin{picture}(60,24)(6,-7)
\put(12,5){\line(1,0){6}}
\put(30,5){\line(1,0){18}}
\put(49,6){\line(1,1){10}}
\put(49,4){\line(1,-1){10}}
\put(21,5){\makebox(0,0){$.$}}
\put(24,5){\makebox(0,0){$.$}}
\put(27,5){\makebox(0,0){$.$}}
\put(12,5){\makebox(0,0){$\bullet$}}
\put(36,5){\makebox(0,0){$\bullet$}}
\put(48,5){\makebox(0,0){$\bullet$}}
\put(60,17){\makebox(0,0){$\bullet$}}
\put(60,-7){\makebox(0,0){$\bullet$}}
\put(12,12){\makebox(0,0){\scriptsize$0$}}
\put(36,12){\makebox(0,0){\scriptsize$0$}}
\put(48,12){\makebox(0,0){\scriptsize$0$}}
\put(65,17){\makebox(0,0){\scriptsize$0$}}
\put(65,-7){\makebox(0,0){\scriptsize$1$}}
\end{picture}}
\]
of ${\mathfrak{so}}(n)$. In any case, the geometric import of these
observations is that there is a complex vector bundle $S$ defined on an
arbitrary spin manifold such that the dif\/ferential operator $D:E\to E$ defined
in~\S\ref{defofD}, when acting on complex-valued sections, becomes
\[
{\mathbb{C}}\otimes_{\mathbb{R}}E\cong S\otimes{\mathbb{C}}^N
\xrightarrow{\,\Dirac\otimes{\mathrm{Id}}\,}S\otimes{\mathbb{C}}^N
\cong{\mathbb{C}}\otimes_{\mathbb{R}}E.
\]
Here ${\mathbb{C}}^N$ denotes the trivial bundle of rank $N=2^{n/2}$ obtained
as the induced bundle from the trivial representation of ${\mathrm{Spin}}(n)$
on ${\mathbb{U}}^*$. In other words, the Dirac operator $D$ from Clif\/ford
analysis may be viewed simply as $2^{n/2}$-copies of an operator
$\Dirac:S\to S$. The operator $\Dirac$ is the classical Dirac operator. There
is, however, an awkward proviso to this conclusion, namely that $D$ and
$\Dirac$ should be acting on complex-valued sections. The theory of real
spin-bundles is quite complicated~\cite{h}. These particular complications are
simply avoided in Clif\/ford analysis by using the real operator $D$ instead.

Similar conclusions hold in the odd-dimensional case:
\[
{\mathbb{CE}}=\Lambda^\bullet{\mathbb{CW}}\cong
\raisebox{-3pt}{\begin{picture}(55,17)(7,0)
\put(12,5){\line(1,0){6}}
\put(30,5){\line(1,0){18}}
\put(49,6){\line(1,0){10}}
\put(49,4){\line(1,0){10}}
\put(21,5){\makebox(0,0){$.$}}
\put(24,5){\makebox(0,0){$.$}}
\put(27,5){\makebox(0,0){$.$}}
\put(12,5){\makebox(0,0){$\bullet$}}
\put(36,5){\makebox(0,0){$\bullet$}}
\put(48,5){\makebox(0,0){$\bullet$}}
\put(60,5){\makebox(0,0){$\bullet$}}
\put(54,5){\makebox(0,0){$\rangle$}}
\put(12,12){\makebox(0,0){\scriptsize$0$}}
\put(36,12){\makebox(0,0){\scriptsize$0$}}
\put(48,12){\makebox(0,0){\scriptsize$0$}}
\put(60,12){\makebox(0,0){\scriptsize$1$}}
\end{picture}}\otimes{\mathbb{C}}^N\enskip\implies\enskip
{\mathbb{C}}\otimes_{\mathbb{R}}E\cong
S\otimes{\mathbb{C}}^N\quad\mbox{and}\quad D=\Dirac\otimes{\mathrm{Id}},
\]
where $\Dirac$ is the classical Dirac operator in odd dimensions
and~$N=2^{(n+1)/2}$.

In fact, all complications are avoided by using the operator~$D$. This
`one-size-fits-all' approach avoids complex bundles and, at the same time,
there is no need to treat even-dimensional and odd-dimensional manifolds
dif\/ferently. The integral formulae for monogenic functions developed
in~\cite{bds}, for example, are obtained in a uniform dimension-free manner.

The main point of this article is to provide a simple formulation and
na\"{\i}ve computational proof of Theorem~\ref{monogenicinvariance} and its
main consequence, the conformal invariance of monogenic functions. But the
conformal invariance of the classical Dirac operator is another consequence.
In arbitrary dimensions, this invariance was first shown by
Kosmann-Schwarzbach~\cite{k-s}.

\section[Rarita-Schwinger operators and their generalisations]{Rarita--Schwinger operators and their generalisations}

Having defined and established the conformal invariance of the Dirac operator
in the context of Clif\/ford analysis, it is reasonably straightforward to do the
same thing for the classical Rarita--Schwinger operator~\cite{bures} and its
symmetric analogues~\cite{bssv}. Throughout this section ${\mathbb{E}}$ is to
be regarded as an ${\mathfrak{so}}(n)$-module under~(\ref{spinaction}).

For the Rarita--Schwinger operator, define a representation ${\mathbb{F}}$ of
${\mathfrak{so}}(n)$ via the exact sequence
\begin{gather}\label{ses}0\to{\mathbb{F}}\to{\mathbb{W}}\otimes{\mathbb{E}}
\xrightarrow{\epsilon}{\mathbb{E}}\to 0,\end{gather}
bearing in mind that $\epsilon$ is a homomorphism of
${\mathfrak{so}}(n)$-modules in accordance with Proposition~\ref{first}.
Lemma~\ref{metricimage} provides a canonical splitting of this sequence.
Specifically, using abstract indices as in the proof of
Proposition~\ref{second}, we define
$\Pi:{\mathbb{W}}\otimes{\mathbb{E}}\to{\mathbb{F}}$ by
\[
\textstyle T_{a\alpha}\stackrel{\Pi}{\longmapsto}
T_{a\alpha}+
\frac{1}{n}\epsilon_{a\alpha}{}^\beta\epsilon^c{}_\beta{}^\gamma T_{c\gamma}.
\]
Finally, to define the Rarita--Schwinger operator via Theorem~\ref{fundamental},
we consider the ${\mathfrak{so}}(n)$-module homomorphism
$\theta:{\mathbb{W}}\otimes{\mathbb{F}}\to{\mathbb{F}}$ given as the
composition
\[
{\mathbb{W}}\otimes{\mathbb{F}}\hookrightarrow
{\mathbb{W}}\otimes{\mathbb{W}}\otimes{\mathbb{E}}\xrightarrow{\tilde\epsilon}
{\mathbb{W}}\otimes{\mathbb{E}}\xrightarrow{\Pi}{\mathbb{F}}
\]
where $\tilde\epsilon(w\otimes v\otimes e)=v\otimes w.e$ or, using abstract
indices,
\[
\textstyle T_{ab\alpha}\stackrel{\tilde\epsilon}{\longmapsto}
\epsilon^b{}_\alpha{}^\beta T_{ba\beta}.
\]
Let us write $\tau$ for the representation of
${\mathfrak{so}}(n)=\Lambda^2{\mathbb{W}}$ on~${\mathbb{F}}$. By construction,
it is the restriction to ${\mathbb{F}}$ of the action
\begin{gather}\label{combinedaction}
\Lambda^2{\mathbb{W}}\otimes{\mathbb{W}}\otimes{\mathbb{E}}\ni X\otimes
w\otimes e\mapsto Xw\otimes e+w\otimes\sigma(X\otimes e)\in
{\mathbb{W}}\otimes{\mathbb{E}}
\end{gather}
on ${\mathbb{W}}\otimes{\mathbb{E}}$.

\begin{proposition}\label{fourth} The composition
\[
{\mathbb{W}}\otimes{\mathbb{F}}\xrightarrow{\iota\otimes{\mathrm{Id}}}
{\mathbb{W}}\otimes\Lambda^2{\mathbb{W}}\otimes{\mathbb{F}}
\xrightarrow{{\mathrm{Id}}\otimes\tau}{\mathbb{W}}\otimes{\mathbb{F}}
\xrightarrow{\theta}{\mathbb{F}}
\]
is equal to
$-\frac{n-1}2\theta:{\mathbb{W}}\otimes{\mathbb{F}}\to{\mathbb{F}}$.
\end{proposition}

\begin{proof} In full, we may expand this composition as
\[
{\mathbb{W}}\otimes{\mathbb{F}}\hookrightarrow
\underbrace{{\mathbb{W}}\otimes{\mathbb{W}}\otimes{\mathbb{E}}
\xrightarrow{\iota\otimes{\mathrm{Id}}\otimes{\mathrm{Id}}}
{\mathbb{W}}\otimes\Lambda^2{\mathbb{W}}\otimes{\mathbb{W}}\otimes{\mathbb{E}}
\xrightarrow{{\mathrm{Id}}\otimes\tau}
{\mathbb{W}}\otimes{\mathbb{W}}\otimes{\mathbb{E}}}_{\mbox{$\star$}}
\xrightarrow{\tilde\epsilon}
{\mathbb{W}}\otimes{\mathbb{E}}\xrightarrow{\Pi}{\mathbb{F}}
\]
and, according to~(\ref{combinedaction}), the homomorphism $\star$ is the sum
of two parts, namely
\begin{gather}\label{one}
{\mathbb{W}}\otimes{\mathbb{W}}\otimes\framebox{${\mathbb{E}}$}
\xrightarrow{\iota\otimes{\mathrm{Id}}\otimes{\mathrm{Id}}}
{\mathbb{W}}\otimes\Lambda^2{\mathbb{W}}\otimes{\mathbb{W}}\otimes
\framebox{${\mathbb{E}}$}
\longrightarrow
{\mathbb{W}}\otimes{\mathbb{W}}\otimes\framebox{${\mathbb{E}}$}\end{gather}
and
\begin{gather}\label{two}
{\mathbb{W}}\otimes\framebox{${\mathbb{W}}$}\otimes{\mathbb{E}}
\xrightarrow{\iota\otimes{\mathrm{Id}}\otimes{\mathrm{Id}}}
{\mathbb{W}}\otimes\Lambda^2{\mathbb{W}}\otimes
\framebox{${\mathbb{W}}$}\otimes{\mathbb{E}}
\longrightarrow
{\mathbb{W}}\otimes\framebox{${\mathbb{W}}$}\otimes{\mathbb{E}},\end{gather}
in which the boxed vector spaces are passengers. The first one is
\[
T_{ab\alpha}\mapsto-2g_{a[b}T_{c]d\alpha}\mapsto
T_{ba\alpha}-g_{ab}T^c{}_{c\alpha}
\]
and composing with $\tilde\epsilon$ gives
\[
T_{ab\alpha}\mapsto
\epsilon^b{}_\alpha{}^\beta(T_{ab\beta}-g_{ab}T^c{}_{c\beta})=
-\epsilon_{a\alpha}{}^\beta T^d{}_{d\beta},
\]
when acting on~${\mathbb{W}}\otimes{\mathbb{F}}$. Evidently, this is in the
kernel of $\Pi$ and so (\ref{one}) makes no contribution to the overall
composition. On the other hand (\ref{two}) can be continued to
\[
{\mathbb{W}}\otimes\framebox{${\mathbb{W}}$}\otimes{\mathbb{E}}
\xrightarrow{\iota\otimes{\mathrm{Id}}\otimes{\mathrm{Id}}}
{\mathbb{W}}\otimes\Lambda^2{\mathbb{W}}\otimes
\framebox{${\mathbb{W}}$}\otimes{\mathbb{E}} \longrightarrow
{\mathbb{W}}\otimes\framebox{${\mathbb{W}}$}\otimes{\mathbb{E}}
\xrightarrow{\tilde\epsilon}
\framebox{${\mathbb{W}}$}\otimes{\mathbb{E}},
\]
which has already been computed in Proposition~\ref{second}. It is
$-\frac{n-1}2\tilde\epsilon$. Therefore, composing with $\Pi$ gives
$-\frac{n-1}2\theta$, as advertised.
\end{proof}

\begin{theorem}
Let $F$ denote the bundle defined on a conformal spin manifold corresponding to
the representation ${\mathbb{F}}$ of\/ ${\mathfrak{so}}(n)$. Then there is a
conformally invariant first order linear differential operator
\[
\textstyle F[-\frac{n-1}2]\xrightarrow{D}F[-\frac{n+1}2]
\]
whose symbol is induced by
$\theta:{\mathbb{W}}\otimes{\mathbb{F}}\to{\mathbb{F}}$.
\end{theorem}
\begin{proof}
Proposition~\ref{fourth} ensures that the criterion of
Theorem~\ref{fundamental} is satisfied.
\end{proof}
The operators in this theorem are the Rarita--Schwinger operators in the context
of Clif\/ford analysis. The symmetric analogues of \cite{bssv} are obtained by a
similar construction starting with the exact sequence
\[
\textstyle 0\to{\mathbb{F}}_j\to\bigodot^j\!{\mathbb{W}}\otimes{\mathbb{E}}
\xrightarrow{{\mathrm{Id}}\otimes\epsilon}
\bigodot^{j-1}\!{\mathbb{W}}\otimes{\mathbb{E}}\to 0,
\]
generalising~(\ref{ses}). Details are left to the reader.

\section{Flat structures and comparison of notations}
The operators that we have constructed in previous sections are conformally
invariant in the `curved' setting. Explicitly, this means that for each
operator there is a universal formula in terms of a chosen Riemannian metric
and its Levi-Civita connection so that, for an arbitrary metric, using the same
formula with any metric in the same conformal class gives the same operator.
This is a very strong notion of invariance. There are several weaker notions
one of which is to restrict attention to conformally f\/lat metrics, asking only
for a local formula in terms of a f\/lat metric from the conformal class and only
that the result be invariant under arbitrary f\/lat-to-f\/lat conformal rescalings.
Evidently, an invariant operator in the fully curved sense gives rise to an
invariant operator in the conformally f\/lat sense as just defined. That
curved invariance is strictly stronger, however, was demonstrated by
Graham~\cite{g} who showed that the operator
\[
\Delta^3:\Lambda^0[1]\to\Lambda^0[-5]\qquad\mbox{on} \quad {\mathbb{R}}^4,
\]
where $\Delta$ is the Laplacian for the standard f\/lat metric on
${\mathbb{R}}^4$ does not arise from a curved invariant operator although it is
easily verified to be invariant under f\/lat-to-f\/lat conformal rescalings (see
also~\cite{gh}).

In any case, the Dirac operator is conformally invariant in the curved case and
therefore invariant under f\/lat-to-f\/lat conformal rescalings. Such rescalings
are scarce. In fact, the only way that they can arise is by so-called `M\"obius
transformations', i.e.\ the transformations obtained from the action of
${\mathrm{SO}}(n+1,1)$ on the $n$-sphere and viewed in ${\mathbb{R}}^n$ by
stereographic projection. Hence, monogenic functions locally defined on
${\mathbb{R}}^n$ are invariant under M\"obius transformations. This is a
well-known phenomenon from Clif\/ford analysis usually attributed to the
invariance of the genera\-lised Cauchy integral formula~\cite{bds} on
${\mathbb{R}}^n$ under M\"obius transformations. Monogenic functions on
conformally f\/lat spin mani\-folds such as cylinders and tori provide a natural
extension~\cite{kr1,kr2} of automorphic forms from the realm of complex
analysis. The $n$-sphere is a spin mani\-fold homogeneous under its conformal
motions, finitely covered by the group $G={\mathrm{SO}}(n+1,1)$ (as indicated below
and fully explained, for example, in~\cite{eg}). It follows that the Dirac
operator is $G$-invariant on~$S^n$ (some authors write $G$-equivariant). This
is yet another sense of `conformal invariance' and a notion that is purely
algebraic~-- the question of classi\-fying the $G$-invariant operators on $S^n$
acting between irreducible homogeneous vector bundles is a~(solved) question in
representation theory. It is this sense of conformal invariance that is in
ef\/fect throughout the following discussion.

Finally, for completeness, we present the conformal invariance of the Dirac
operator on the f\/lat model in terms of the Dynkin diagram notation
of~\cite{beastwood}. The conformal sphere $S^n$ is a~homogeneous space for the
Lie group ${\mathrm{SO}}(n+1,1)$. For convenience we shall take this to mean
the connected component of the group of $(n+2)\times(n+2)$ matrices preserving
the quadratic form
\[
\left[\begin{array}{ccc} 0&0&1\\
0&\mbox{${\mathrm{Id}}$}&0\\ 1&0&0
\end{array}\right]
\]
where ${\mathrm{Id}}$ denotes the $n\times n$ identity matrix. This group acts
on $S^n$ as its conformal motions and, for a suitable choice of basepoint,
$S^n={\mathrm{SO}}(n+1,1)/P$ where $P$ is the subgroup consisting of matrices
of the form
\[
\left[\begin{array}{ccc} \lambda&\cdot&\cdot\\
0&\mbox{${\mathrm{M}}$}&\cdot\\ 0&0&\lambda^{-1}
\end{array}\right],\qquad\mbox{for}\quad
\left\lbrace\begin{array}l\lambda>0\\ M\in{\mathrm{SO}}(n).\end{array}\right.
\]
The irreducible real homogeneous bundles on $S^n$ are thus parameterised by the
representations
\[
\left[\begin{array}{ccc} \lambda&\cdot&\cdot\\
0&\mbox{${\mathrm{M}}$}&\cdot\\ 0&0&\lambda^{-1}
\end{array}\right]\mapsto\lambda^{-w}\rho(M)
\]
for $w\in{\mathbb{R}}$ and $\rho$ an irreducible real representation
of~${\mathrm{SO}}(n)$. It is convenient to specify such a~representation of
$P$ by recording its highest (or lowest) weight as a linear combination of the
fundamental weights for ${\mathfrak{so}}(n+1,1)$. With the conventions of
\cite{beastwood} we find, for example, that
\begin{gather*}
\hspace*{1mm}\Lambda^0[w]=\raisebox{-3pt}{\begin{picture}(51,17)(-3,0)
\put(0,5){\line(1,0){36}}
\put(37,6){\line(1,0){10}}
\put(37,4){\line(1,0){10}}
\put(0,5){\makebox(0,0){$\times$}}
\put(12,5){\makebox(0,0){$\bullet$}}
\put(24,5){\makebox(0,0){$\bullet$}}
\put(36,5){\makebox(0,0){$\bullet$}}
\put(48,5){\makebox(0,0){$\bullet$}}
\put(42,5){\makebox(0,0){$\rangle$}}
\put(0,12){\makebox(0,0){\scriptsize$w$}}
\put(12,12){\makebox(0,0){\scriptsize$1$}}
\put(24,12){\makebox(0,0){\scriptsize$0$}}
\put(36,12){\makebox(0,0){\scriptsize$0$}}
\put(48,12){\makebox(0,0){\scriptsize$0$}}
\end{picture}}
\qquad
\Lambda^1=\raisebox{-3pt}{\begin{picture}(51,17)(-3,0)
\put(0,5){\line(1,0){36}}
\put(37,6){\line(1,0){10}}
\put(37,4){\line(1,0){10}}
\put(0,5){\makebox(0,0){$\times$}}
\put(12,5){\makebox(0,0){$\bullet$}}
\put(24,5){\makebox(0,0){$\bullet$}}
\put(36,5){\makebox(0,0){$\bullet$}}
\put(48,5){\makebox(0,0){$\bullet$}}
\put(42,5){\makebox(0,0){$\rangle$}}
\put(-1,12){\makebox(0,0){\scriptsize$-2$}}
\put(12,12){\makebox(0,0){\scriptsize$1$}}
\put(24,12){\makebox(0,0){\scriptsize$0$}}
\put(36,12){\makebox(0,0){\scriptsize$0$}}
\put(48,12){\makebox(0,0){\scriptsize$0$}}
\end{picture}}
\qquad
\Lambda^2=\raisebox{-3pt}{\begin{picture}(51,17)(-3,0)
\put(0,5){\line(1,0){36}}
\put(37,6){\line(1,0){10}}
\put(37,4){\line(1,0){10}}
\put(0,5){\makebox(0,0){$\times$}}
\put(12,5){\makebox(0,0){$\bullet$}}
\put(24,5){\makebox(0,0){$\bullet$}}
\put(36,5){\makebox(0,0){$\bullet$}}
\put(48,5){\makebox(0,0){$\bullet$}}
\put(42,5){\makebox(0,0){$\rangle$}}
\put(-1,12){\makebox(0,0){\scriptsize$-3$}}
\put(12,12){\makebox(0,0){\scriptsize$0$}}
\put(24,12){\makebox(0,0){\scriptsize$1$}}
\put(36,12){\makebox(0,0){\scriptsize$0$}}
\put(48,12){\makebox(0,0){\scriptsize$0$}}
\end{picture}}
\qquad
\Lambda^1[w]=\raisebox{-3pt}{\begin{picture}(59,17)(-11,0)
\put(-4,5){\line(1,0){40}}
\put(37,6){\line(1,0){10}}
\put(37,4){\line(1,0){10}}
\put(-4,5){\makebox(0,0){$\times$}}
\put(12,5){\makebox(0,0){$\bullet$}}
\put(24,5){\makebox(0,0){$\bullet$}}
\put(36,5){\makebox(0,0){$\bullet$}}
\put(48,5){\makebox(0,0){$\bullet$}}
\put(42,5){\makebox(0,0){$\rangle$}}
\put(-4,12){\makebox(0,0){\scriptsize$w\!-\!2$}}
\put(12,12){\makebox(0,0){\scriptsize$1$}}
\put(24,12){\makebox(0,0){\scriptsize$0$}}
\put(36,12){\makebox(0,0){\scriptsize$0$}}
\put(48,12){\makebox(0,0){\scriptsize$0$}}
\end{picture}}
\end{gather*}
on~$S^9$ and
\begin{gather*}
\hspace*{1mm}\Lambda^0[w]=\raisebox{-10pt}{\begin{picture}(62,24)(-2,-7)
\put(0,5){\line(1,0){36}}
\put(37,6){\line(1,1){10}}
\put(37,4){\line(1,-1){10}}
\put(0,5){\makebox(0,0){$\times$}}
\put(12,5){\makebox(0,0){$\bullet$}}
\put(24,5){\makebox(0,0){$\bullet$}}
\put(36,5){\makebox(0,0){$\bullet$}}
\put(48,17){\makebox(0,0){$\bullet$}}
\put(48,-7){\makebox(0,0){$\bullet$}}
\put(0,12){\makebox(0,0){\scriptsize$w$}}
\put(12,12){\makebox(0,0){\scriptsize$1$}}
\put(24,12){\makebox(0,0){\scriptsize$0$}}
\put(36,12){\makebox(0,0){\scriptsize$0$}}
\put(53,17){\makebox(0,0){\scriptsize$0$}}
\put(53,-7){\makebox(0,0){\scriptsize$0$}}
\end{picture}}
\quad
\Lambda^1=\raisebox{-10pt}{\begin{picture}(62,24)(-2,-7)
\put(0,5){\line(1,0){36}}
\put(37,6){\line(1,1){10}}
\put(37,4){\line(1,-1){10}}
\put(0,5){\makebox(0,0){$\times$}}
\put(12,5){\makebox(0,0){$\bullet$}}
\put(24,5){\makebox(0,0){$\bullet$}}
\put(36,5){\makebox(0,0){$\bullet$}}
\put(48,17){\makebox(0,0){$\bullet$}}
\put(48,-7){\makebox(0,0){$\bullet$}}
\put(-1,12){\makebox(0,0){\scriptsize$-2$}}
\put(12,12){\makebox(0,0){\scriptsize$1$}}
\put(24,12){\makebox(0,0){\scriptsize$0$}}
\put(36,12){\makebox(0,0){\scriptsize$0$}}
\put(53,17){\makebox(0,0){\scriptsize$0$}}
\put(53,-7){\makebox(0,0){\scriptsize$0$}}
\end{picture}}
\quad
\Lambda^2=\raisebox{-10pt}{\begin{picture}(62,24)(-2,-7)
\put(0,5){\line(1,0){36}}
\put(37,6){\line(1,1){10}}
\put(37,4){\line(1,-1){10}}
\put(0,5){\makebox(0,0){$\times$}}
\put(12,5){\makebox(0,0){$\bullet$}}
\put(24,5){\makebox(0,0){$\bullet$}}
\put(36,5){\makebox(0,0){$\bullet$}}
\put(48,17){\makebox(0,0){$\bullet$}}
\put(48,-7){\makebox(0,0){$\bullet$}}
\put(-1,12){\makebox(0,0){\scriptsize$-3$}}
\put(12,12){\makebox(0,0){\scriptsize$0$}}
\put(24,12){\makebox(0,0){\scriptsize$1$}}
\put(36,12){\makebox(0,0){\scriptsize$0$}}
\put(53,17){\makebox(0,0){\scriptsize$0$}}
\put(53,-7){\makebox(0,0){\scriptsize$0$}}
\end{picture}}
\quad
\Lambda^1[w]=\raisebox{-10pt}{\begin{picture}(62,24)(-10,-7)
\put(-4,5){\line(1,0){40}}
\put(37,6){\line(1,1){10}}
\put(37,4){\line(1,-1){10}}
\put(-4,5){\makebox(0,0){$\times$}}
\put(12,5){\makebox(0,0){$\bullet$}}
\put(24,5){\makebox(0,0){$\bullet$}}
\put(36,5){\makebox(0,0){$\bullet$}}
\put(48,17){\makebox(0,0){$\bullet$}}
\put(48,-7){\makebox(0,0){$\bullet$}}
\put(-4,12){\makebox(0,0){\scriptsize$w\!-\!2$}}
\put(12,12){\makebox(0,0){\scriptsize$1$}}
\put(24,12){\makebox(0,0){\scriptsize$0$}}
\put(36,12){\makebox(0,0){\scriptsize$0$}}
\put(53,17){\makebox(0,0){\scriptsize$0$}}
\put(53,-7){\makebox(0,0){\scriptsize$0$}}
\end{picture}}
\end{gather*}
on~$S^{10}$. On odd-dimensional spheres $S^n$ the Dirac operator acts as
\begin{gather*}
\hspace*{5mm}
\raisebox{-3pt}{\begin{picture}(65,17)(-3,0)
\put(0,5){\line(1,0){18}}
\put(30,5){\line(1,0){18}}
\put(49,6){\line(1,0){10}}
\put(49,4){\line(1,0){10}}
\put(21,5){\makebox(0,0){$.$}}
\put(24,5){\makebox(0,0){$.$}}
\put(27,5){\makebox(0,0){$.$}}
\put(0,5){\makebox(0,0){$\times$}}
\put(12,5){\makebox(0,0){$\bullet$}}
\put(36,5){\makebox(0,0){$\bullet$}}
\put(48,5){\makebox(0,0){$\bullet$}}
\put(60,5){\makebox(0,0){$\bullet$}}
\put(54,5){\makebox(0,0){$\rangle$}}
\put(-4,12){\makebox(0,0){\scriptsize$-n/2$}}
\put(12,12){\makebox(0,0){\scriptsize$0$}}
\put(36,12){\makebox(0,0){\scriptsize$0$}}
\put(48,12){\makebox(0,0){\scriptsize$0$}}
\put(60,12){\makebox(0,0){\scriptsize$1$}}
\end{picture}}\longrightarrow
\raisebox{-3pt}{\begin{picture}(82,17)(-18,0)
\put(-10,5){\line(1,0){28}}
\put(30,5){\line(1,0){18}}
\put(49,6){\line(1,0){10}}
\put(49,4){\line(1,0){10}}
\put(21,5){\makebox(0,0){$.$}}
\put(24,5){\makebox(0,0){$.$}}
\put(27,5){\makebox(0,0){$.$}}
\put(-10,5){\makebox(0,0){$\times$}}
\put(12,5){\makebox(0,0){$\bullet$}}
\put(36,5){\makebox(0,0){$\bullet$}}
\put(48,5){\makebox(0,0){$\bullet$}}
\put(60,5){\makebox(0,0){$\bullet$}}
\put(54,5){\makebox(0,0){$\rangle$}}
\put(-10,12){\makebox(0,0){\scriptsize$-n/2\!-\!1$}}
\put(12,12){\makebox(0,0){\scriptsize$0$}}
\put(36,12){\makebox(0,0){\scriptsize$0$}}
\put(48,12){\makebox(0,0){\scriptsize$0$}}
\put(60,12){\makebox(0,0){\scriptsize$1$}}
\end{picture}}.
\end{gather*}
On even-dimensional spheres, however, there are two irreducible Dirac
operators:
\begin{gather*}
\hspace*{5mm}
\raisebox{-10pt}{\begin{picture}(68,24)(-2,-7)
\put(0,5){\line(1,0){18}}
\put(30,5){\line(1,0){18}}
\put(49,6){\line(1,1){10}}
\put(49,4){\line(1,-1){10}}
\put(21,5){\makebox(0,0){$.$}}
\put(24,5){\makebox(0,0){$.$}}
\put(27,5){\makebox(0,0){$.$}}
\put(0,5){\makebox(0,0){$\times$}}
\put(12,5){\makebox(0,0){$\bullet$}}
\put(36,5){\makebox(0,0){$\bullet$}}
\put(48,5){\makebox(0,0){$\bullet$}}
\put(60,17){\makebox(0,0){$\bullet$}}
\put(60,-7){\makebox(0,0){$\bullet$}}
\put(-4,12){\makebox(0,0){\scriptsize$-n/2$}}
\put(12,12){\makebox(0,0){\scriptsize$0$}}
\put(36,12){\makebox(0,0){\scriptsize$0$}}
\put(48,12){\makebox(0,0){\scriptsize$0$}}
\put(65,17){\makebox(0,0){\scriptsize$1$}}
\put(65,-7){\makebox(0,0){\scriptsize$0$}}
\end{picture}}
\longrightarrow
\raisebox{-10pt}{\begin{picture}(83,24)(-17,-7)
\put(-10,5){\line(1,0){28}}
\put(30,5){\line(1,0){18}}
\put(49,6){\line(1,1){10}}
\put(49,4){\line(1,-1){10}}
\put(21,5){\makebox(0,0){$.$}}
\put(24,5){\makebox(0,0){$.$}}
\put(27,5){\makebox(0,0){$.$}}
\put(-10,5){\makebox(0,0){$\times$}}
\put(12,5){\makebox(0,0){$\bullet$}}
\put(36,5){\makebox(0,0){$\bullet$}}
\put(48,5){\makebox(0,0){$\bullet$}}
\put(60,17){\makebox(0,0){$\bullet$}}
\put(60,-7){\makebox(0,0){$\bullet$}}
\put(-10,12){\makebox(0,0){\scriptsize$-n/2\!-\!1$}}
\put(12,12){\makebox(0,0){\scriptsize$0$}}
\put(36,12){\makebox(0,0){\scriptsize$0$}}
\put(48,12){\makebox(0,0){\scriptsize$0$}}
\put(65,17){\makebox(0,0){\scriptsize$0$}}
\put(65,-7){\makebox(0,0){\scriptsize$1$}}
\end{picture}}
\quad\mbox{and}\qquad
\raisebox{-10pt}{\begin{picture}(68,24)(-2,-7)
\put(0,5){\line(1,0){18}}
\put(30,5){\line(1,0){18}}
\put(49,6){\line(1,1){10}}
\put(49,4){\line(1,-1){10}}
\put(21,5){\makebox(0,0){$.$}}
\put(24,5){\makebox(0,0){$.$}}
\put(27,5){\makebox(0,0){$.$}}
\put(0,5){\makebox(0,0){$\times$}}
\put(12,5){\makebox(0,0){$\bullet$}}
\put(36,5){\makebox(0,0){$\bullet$}}
\put(48,5){\makebox(0,0){$\bullet$}}
\put(60,17){\makebox(0,0){$\bullet$}}
\put(60,-7){\makebox(0,0){$\bullet$}}
\put(-4,12){\makebox(0,0){\scriptsize$-n/2$}}
\put(12,12){\makebox(0,0){\scriptsize$0$}}
\put(36,12){\makebox(0,0){\scriptsize$0$}}
\put(48,12){\makebox(0,0){\scriptsize$0$}}
\put(65,17){\makebox(0,0){\scriptsize$0$}}
\put(65,-7){\makebox(0,0){\scriptsize$1$}}
\end{picture}}
\longrightarrow
\raisebox{-10pt}{\begin{picture}(87,24)(-17,-7)
\put(-10,5){\line(1,0){28}}
\put(30,5){\line(1,0){18}}
\put(49,6){\line(1,1){10}}
\put(49,4){\line(1,-1){10}}
\put(21,5){\makebox(0,0){$.$}}
\put(24,5){\makebox(0,0){$.$}}
\put(27,5){\makebox(0,0){$.$}}
\put(-10,5){\makebox(0,0){$\times$}}
\put(12,5){\makebox(0,0){$\bullet$}}
\put(36,5){\makebox(0,0){$\bullet$}}
\put(48,5){\makebox(0,0){$\bullet$}}
\put(60,17){\makebox(0,0){$\bullet$}}
\put(60,-7){\makebox(0,0){$\bullet$}}
\put(-10,12){\makebox(0,0){\scriptsize$-n/2\!-\!1$}}
\put(12,12){\makebox(0,0){\scriptsize$0$}}
\put(36,12){\makebox(0,0){\scriptsize$0$}}
\put(48,12){\makebox(0,0){\scriptsize$0$}}
\put(65,17){\makebox(0,0){\scriptsize$1$}}
\put(65,-7){\makebox(0,0){\scriptsize$0$}}
\end{picture}}.
\end{gather*}
In odd dimensions the generalised Rarita--Schwinger operators of \cite{bssv} act
as
\begin{gather*}
\hspace*{15mm}
\raisebox{-3pt}{\begin{picture}(65,17)(-3,0)
\put(-24,5){\line(1,0){42}}
\put(30,5){\line(1,0){18}}
\put(49,6){\line(1,0){10}}
\put(49,4){\line(1,0){10}}
\put(21,5){\makebox(0,0){$.$}}
\put(24,5){\makebox(0,0){$.$}}
\put(27,5){\makebox(0,0){$.$}}
\put(-24,5){\makebox(0,0){$\times$}}
\put(0,5){\makebox(0,0){$\bullet$}}
\put(12,5){\makebox(0,0){$\bullet$}}
\put(36,5){\makebox(0,0){$\bullet$}}
\put(48,5){\makebox(0,0){$\bullet$}}
\put(60,5){\makebox(0,0){$\bullet$}}
\put(54,5){\makebox(0,0){$\rangle$}}
\put(-26,12){\makebox(0,0){\scriptsize$-n/2\!-\!j$}}
\put(0,12){\makebox(0,0){\scriptsize$j$}}
\put(12,12){\makebox(0,0){\scriptsize$0$}}
\put(36,12){\makebox(0,0){\scriptsize$0$}}
\put(48,12){\makebox(0,0){\scriptsize$0$}}
\put(60,12){\makebox(0,0){\scriptsize$1$}}
\end{picture}}\longrightarrow
\raisebox{-3pt}{\begin{picture}(108,17)(-44,0)
\put(-30,5){\line(1,0){48}}
\put(30,5){\line(1,0){18}}
\put(49,6){\line(1,0){10}}
\put(49,4){\line(1,0){10}}
\put(21,5){\makebox(0,0){$.$}}
\put(24,5){\makebox(0,0){$.$}}
\put(27,5){\makebox(0,0){$.$}}
\put(-30,5){\makebox(0,0){$\times$}}
\put(0,5){\makebox(0,0){$\bullet$}}
\put(12,5){\makebox(0,0){$\bullet$}}
\put(36,5){\makebox(0,0){$\bullet$}}
\put(48,5){\makebox(0,0){$\bullet$}}
\put(60,5){\makebox(0,0){$\bullet$}}
\put(54,5){\makebox(0,0){$\rangle$}}
\put(-30,12){\makebox(0,0){\scriptsize$-n/2\!-\!j\!-\!1$}}
\put(0,12){\makebox(0,0){\scriptsize$j$}}
\put(12,12){\makebox(0,0){\scriptsize$0$}}
\put(36,12){\makebox(0,0){\scriptsize$0$}}
\put(48,12){\makebox(0,0){\scriptsize$0$}}
\put(60,12){\makebox(0,0){\scriptsize$1$}}
\end{picture}}.
\end{gather*}

\subsection*{Acknowledgements}
It is a pleasure to acknowledgment useful conversations with Vladim\'{\i}r
Sou\v{c}ek. He is certainly one person for whom the `well-known' material in
this article is actually known.

Michael Eastwood is a Professorial Fellow of the Australian Research Council.
This research was begun during a visit by John Ryan to the University of
Adelaide in 2005, which was also supported by the Australian Research Council.
This support is gratefully acknowledged. John Ryan also thanks the University
of Adelaide for hospitality during his visit.

\pdfbookmark[1]{References}{ref}
\LastPageEnding

\end{document}